\documentclass[12pt,a4paper]{article}

\usepackage{amsmath}
\usepackage{amssymb, latexsym}
\usepackage{color}
\newcommand{\ourEE}{\mathbb E}
\newcommand{\nn}{\textbf n}

\newcommand{\ourN}{\mbox{$\mathbb N$}}
\newcommand{{\Z}}{{\bf Z}}
\newcommand{\R}{\mbox{$\mathbb R$}}

\newcommand{\omfp}{(\mbox{$\Omega$,$\cal F,$$P)$}}

\newcommand{\ourE}{\mathbb E}

\newtheorem{theorem}{Theorem}
\newtheorem{cor}[theorem]{Corollary}
\newtheorem{defn}[theorem]{Definition}
\newtheorem{lemma}[theorem]{Lemma}
\newtheorem{prop}[theorem]{Proposition}
\newtheorem{ex}[theorem]{Example}
\newtheorem{rem}[theorem]{Remark}

\newtheorem{assumption}[theorem]{Assumption}

\numberwithin{theorem}{section}
\numberwithin{equation}{section}

\newcommand{\beq}{\begin{equation}}
\newcommand{\eeq}{\end{equation}}
\newcommand{\beqo}{\begin{equation*}}
\newcommand{\eeqo}{\end{equation*}}
\newcommand{\bdm}{\begin{displaymath}}
\newcommand{\edm}{\end{displaymath}}
\newcommand{\beqar}{\begin{eqnarray}}
\newcommand{\eeqar}{\end{eqnarray}}
\newcommand{\beqaro}{\begin{eqnarray*}}
\newcommand{\eeqaro}{\end{eqnarray*}}
\newcommand{\bal}{\begin{align}}
\newcommand{\eal}{\end{align}}

\newcommand{\btheorem}{\begin{theorem}}
\newcommand{\etheorem}{\end{theorem}}
\newcommand{\blemma}{\begin{lemma}}
\newcommand{\elemma}{\end{lemma}}
\newcommand{\bremark}{\begin{rem}}
\newcommand{\eremark}{\end{rem}}
\newcommand{\bcor}{\begin{cor}}
\newcommand{\ecor}{\end{cor}}
\newcommand{\bex}{\begin{ex}}
\newcommand{\eex}{\end{ex}}
\newcommand{\bdefn}{\begin{defn}}
\newcommand{\edefn}{\end{defn}}
\newcommand{\bprop}{\begin{prop}}
\newcommand{\eprop}{\end{prop}}
\newcommand{\bproof}{Proof.\ }
\newcommand{\eproof}{$\Box$\\}
\newcommand{\bass}{\begin{assumption}}
\newcommand{\eass}{\end{assumption}}

\newcommand{\cE}{{\cal E}}
\newcommand{\cF}{{\cal F}}

\newcommand{\n}{\underline n}
\newcommand{\m}{\underline m}

\newcommand{\X}{{\bf X}}
\newcommand{\Y}{{\bf Y}}
\begin{document}

\title{Correlation formulas for Markovian  network processes in a random environment}
\author{Hans Daduna
\thanks{Work supported by Deutscher Akademischer Austauschdienst}
\\
    {\it  Hamburg University}\\
         Department of Mathematics\\
          Bundesstrasse 55\\
          20146 Hamburg\\
          Germany\\
\and Ryszard Szekli
  \thanks{ Work supported by NCN Research Grant DEC-2011/01/B/ST1/01305
 }
 \\
  {\it Wroc{\l}aw University}\\
        Mathematical Institute\\
          pl. Grunwaldzki 2/4\\
          50--384 Wroc{\l}aw\\
           Poland\\
}
\maketitle

  \footnotetext{
  {\it AMS (1991) subject classification}:  60K25, 60J25\\
 {\it Short title}: correlations in  networks\\
 corr-unrel-networks-171013.tex
  }





\begin{abstract}
We consider  Markov  processes, which describe e.g. queueing network processes,
in a random environment which influences the network by determining random
breakdown of nodes, and the necessity of repair thereafter.
Starting from an explicit steady state distribution of product form available
in the literature, we notice that this steady state distribution does not provide
information about the correlation structure in time and space (over nodes). We  study this correlation structure via one step correlations for the
queueing-environment process. Although formulas for absolute values of these  correlations are  complicated, the differences of correlations of related
networks are  simple and have a nice structure.
We therefore compare two networks in a random environment  having the same invariant
distribution, and focus on  the time behaviour of the processes
 when in such a network the environment changes or the
rules for traveling  are perturbed.
Evaluating the comparison formulas we  compare spectral gaps
and asymptotic variances of  related processes.
\end{abstract}



\section{Introduction}\label{sect:intro} 
We consider   stochastic networks of the Jacksonian type in a random environment.
For a general introduction to  Markov processes in  random environments
with applications to  networks, see \cite{zhu:94}, \cite{economou:05},
\cite{balsamo;marin:13}.
These stochastic networks have recently found interest as a general model for queueing
 networks in connection with other areas of Operations Research, e.g. inventory theory and reliability theory. The interaction of network and environment is  (i) in the queueing-inventory
system that the service process  decreases inventories and
the inventory restricts  serving customers due to limited stock at hand, and
(ii) in the queueing-reliability system that  external forces let servers
break down which requires repair.\\
We concentrate in the present paper on the second framework: External forces
 generate random breakdowns of servers in the network and the subsequent repair
is  also  performed  under random influences.
We allow the environment to be of a rather general structure, which implies that nodes may break down in isolation or in groups, and that batch repair is possible as well.\\
For this framework there is a product form extension of Jackson's steady state result
at hand, which provides in case of ergodicity the joint steady state distribution of
the environment (represented by the set of broken down nodes) and the joint queue
length vector in a product form: The environment status and the queue lengths
seem to decouple asymptotically and in steady state (which is the essence of Jackson's theorem
in case of pure queueing systems).\\
Clearly, this does not mean that  environment and  queue lengths are independent:
The environment is assumed here to be a Markov process for its own, but it strongly
influences the service provided by the nodes and even the arrival streams there and,
furthermore, the nodes interact as well - the interactions are carried by the traveling
customers.
These dependencies are not expressed by the one-dimensional (in time) marginal process
distribution, which is a multidimensional (in space) product form distribution.
In fact, very little is known about the dependence structure of the interacting processes.
Therefore we study the correlation structure in time of the
environment-queue length process via the one-step correlations (the Dirichlet forms of the associated Markov process), which in time as well as in
space exhibit complex dependence behaviour.\\
To be more precise, our  interest is focused on the following scenarios:
Compare two networks in a random environment which have the same invariant product form distribution, and so are in some sense variants of one another. Typical questions are: What happens to the time behaviour of a network when in such a network the
rules for traveling (routing chains) are perturbed, or, when the environment changes?\\
Our main results are comparison theorems and formulas which provide differences of
one step correlations
in related, resp. perturbed networks. Although the formulas for absolute values of the one
step correlations are rather complicated, it turns out that differences of correlations of related
networks are surprisingly simple and have a nice structure. As a consequence, whenever we have obtained quantities connected
to one step correlations for some reference network as an anchor
(possibly from simulations or numerical evaluations), we
can perform easily explicit performance analysis, especially sensitivity analysis by
varying, e.g., breakdown and repair probabilities or routing probabilities.

The structure of the paper is as follows. In Section \ref{sect:networks} we describe 
stochastic networks and the influence of the environment via Markovian breakdown and
repair processes, which results in a non-Markovian structure of the queue size processes alone,
and cite the steady state distribution for these networks.
In Section \ref{sect:1stepcorr} we derive the explicit formulas for the one step correlations in time
for the joint environment-network process and show that for the interesting
comparison problems these formulas simplify considerably.
In Section \ref{sect:applications} we show that our results allow to compare the spectral
gaps and asymptotic variances of different systems by evaluating our previous formulas
suitably. Comparison results for spectral gaps allow to compare speed of convergence to stationarity for networks in $L^2$ norm. 
In Section \ref{sect:finite-buffers} we discuss relations  to networks with finite buffers
and extensions of our theorems to this area.
Section \ref{sect:proofs} comprises the main technical proofs.\\
{\bf Notation and conventions:}\\
For a set $M$ we denote by $2^M = {\cal P}(M)$ the set of all subsets of $M$.\\
For sets $A,B$ we write $A\subseteq B$ for $A$ which is a  subset of $B$ or equals $B$, and
we write $A\subset B$ for $A$ which  is a  subset of $B$ but does not equal $B$.\\
$\delta_{xy}$ is the Kronecker Delta, which is $1$ iff $x=y$ and $0$ otherwise.\\
Throughout, the node set of our graphs (networks) are denoted by $\tilde J :=\{1,\dots,J\}$,
and the ''extended node set'' is $\tilde J_0 := \{0,1,\dots,J\}$, where ''$0$'' refers to the
external source and sink of the network.\\
We denote the diagonal matrix with a vector $\xi$ on the diagonal and zero otherwise by $diag(\xi)$.\\
$e_j$ is the standard $j$-th base vector in $\ourN^J$ if $1\leq j\leq J$ and $e_0$ is the $J$-dimensional zero vector.\\
$Id$ is the identity matrix of appropriate dimension, defined by the context.\\
For $D\subseteq \tilde J$ and $\n=(n_j:j\in\tilde J)\in \ourN^J$ we write
$\n_D :=(n_j:j\in D)\in \ourN^{|D|}$ and
$\n_{\tilde J\setminus D} :=(n_j:j\in \tilde J\setminus D)\in \ourN^{|\tilde J\setminus D|}$,
and will, as usual, identify
$\n = ((n_j:j\in D),(n_j:j\in \tilde J\setminus D))= (n_j:j\in D,n_j:j\in \tilde J\setminus D)$.\\
Similarly we use for $\ourN^J$-valued random variables with $X_t = X(t) =\n=(n_j:j\in\tilde J)$
the self explaining abbreviations $X_D(t)=(n_j:j\in D)\in \ourN^{|D|}$
and $X(t)=(X_D(t),X_{\tilde J\setminus D}(t))$.\\
For a probability space $(\tilde{\mathbb E}, {\cal E}, \tilde\pi)$  and functions
$f,g:(\tilde{\mathbb E}, {\cal E})\to (\mathbb R, \mathbb B)$
we define the inner product of $f,g$ with respect to $\tilde\pi$,
whenever the following integral exists:
$
\langle f,g\rangle_{\tilde\pi} = \int_{\tilde{\mathbb E}}f(x)\cdot g(x)\ \tilde\pi(dx).
$\\
$L^2:=L^2(\tilde\ourE,\tilde\pi)$ is the space of square integrable functions with
respect to $\tilde\pi$, and $\|f\|_{\tilde\pi}={(\langle f,f\rangle_{\tilde\pi})}^{1/2}$.\\
All random variables occurring  are defined on a common underlying probability space
$(\Omega, \cF, P)$.

\section{Stochastic networks in a random environment}
\label{sect:networks}

\subsection{Stochastic networks}
\label{sect:networksJackson}
A Jackson network \cite{jackson:57} consists of $J$
nodes numbered $1,\dots,J$, where indistinguishable customers arrive, are served, possibly at several
stations, and eventually depart from the network.
The nodes are exponential single servers with state
dependent service rates and with an infinite waiting room under
first--come--first--served (FCFS)
regime. If at node $j$ there are $n_j>0$ customers present, either
in service or waiting, then service is provided there at rate
$\mu_j(n_j)>0$;  we assume $\sup\{\mu_j(k):j\in\{1,\dots,J\},k\in\mathbb N\}< \infty$
and set $\mu_j(0):= 0.$
All  customers follow the same rules.
We shall need later on a slight extension of the standard Jackson network models.
This is described in terms of an irreducible stochastic routing matrix
\beq\label{eq:routing-rejection1}
R=[r_{ij}]_{i,j=0,\ldots,J},
\eeq
 where  the artificial ''node $0$'' represents the source and the sink of all customers.
Strict inequality may hold for $r_{00}\geq 0$, which means that some
arriving customers may be rejected.\\
Customer  arrive in a Poisson stream of intensity $\lambda>0$ which is split
(independently) according to the first row $r_0 := (r_{0i}:i = 0,1,\dots,J)$ of $R$.
Then at nodes $j=1,2,\dots, J$ we observe independent Poisson-$\lambda_j$ arrival streams with
$\lambda_j = \lambda r_{0j}$, while a portion $\lambda_0 = \lambda r_{00}$ of the arriving customers
is rejected (lost).\\
Routing is Markovian,
a customer departing from node $i$ immediately proceeds to node
$j$ with probability $r_{ij}\geq 0$, and departs from the network
with probability $r_{j0}.$\\
Then the {\bf  traffic equations for the admitted customers}
\begin{equation} \label{eq:trafficequation-adm}
\eta_j=\lambda_j + \sum_{i=1}^J \eta_i r_{ij},\quad\quad
j=1,\dots,J,
\end{equation}
have a unique solution  $(\eta_j:j=1,\dots,J).$ Note, that \eqref{eq:trafficequation-adm} only counts for
the {admitted} customers because of  $\lambda_j = \lambda\cdot r_{0j}, j=1,\dots,J,$
and  $\lambda_1+\ldots +\lambda_J = \lambda(1-r_{00})$. In case of  $r_{00}=0$,
$R$ is the so-called extended routing matrix of standard Jackson networks,
see \cite{daduna;szekli:08}[(3.2)].\\
Let $\X=(X_t:t\geq 0)$ denote the vector process recording the
joint queue lengths in the network at time  $t$.
$X_t=(X_1(t),\dots,X_J(t))\in \ourN^J$ reads: at time $t$ there are
$X_j(t)$ customers present at node $j$, either in service or
waiting. The assumptions put on the system imply that $\X$ is a
strong Markov process on state space $\ourN^J$  with generator $Q^\X = (q^\X(\n,\m):\m,\n\in \ourN^J)$
which is given for $g:\ourN^J\to \R$ by
\begin{eqnarray}\label{eq:generatorJ}
(Q^\X g)(\n)&=&
 \sum_{j=1}^J \lambda_j (g(\n + e_j) - g(\n))
+\sum_{j=1}^J \mu_j(n_j) r_{j0} (g(\n - e_j)- g(\n))\nonumber\\
&&+ \sum_{j=1}^J \mu_j(n_j)
\sum_{{i=1}
}^J
r_{ji} (g(\n - e_j+ e_i) - g(\n))
 \end{eqnarray}
$Q^\X $ is a bounded operator because of $\inf_{\n\in \ourN^J} q^\X (\n,\n) > -\infty$.
We assume throughout that $\X$ is ergodic.

For an ergodic network process $\X$ Jackson's theorem \cite{jackson:57} states  that
the unique steady--state and
limiting distribution $\pi$ on $ \ourN^J$ is with normalizing  constants $C(j)$ for marginal distributions of $\X$
\begin{equation}                        \label{eq:steadystateJ}
\pi(\n)=\pi(n_1,\dots,n_J) =
\prod_{j=1}^{J}\left( C(j)^{-1} \prod_{k=1}^{n_j} \frac{\eta_j}{\mu_j(k)}\right)\,,
\quad \n = (n_1,\dots,n_J)\in \ourN^J,
\end{equation}

\subsection{Breakdown-repair processes}
\label{sect:breakdown-repair}
We are interested in stochastic networks, where nodes due to external environment influences
can breakdown and are repaired periodically.
A common situation is that the breakdown-repair process is Markov of its own, and the network reacts on the perturbations driven by the random environment.
To describe these Markovian processes
we consider a set of $J$ stations or devices  (nodes) numbered $1,\dots,J$. Stations are unreliable, break down randomly and are repaired thereafter.  Repair time is random as well.
We assume that the availability status of the system can be described by a homogeneous
Markov process\\
$
\Y=(Y(t):t\geq 0),\quad Y(t) :\omfp\to (2^{\tilde J},{\cal P}({2^{\tilde J}})).
$\\
$Y(t) = D$ indicates that at time $t\geq 0$
the stations  in $D\subseteq \tilde J$ are  down and under repair, while stations in
$\tilde J \setminus D \subseteq \tilde J$ are functioning (''are up'').
The transition rates (breakdown and repair intensities) of  $\Y$ are
\begin{defn}\label{defn:queueindependentbreakdown}
Take any pair of  functions
$
A: 2^{\tilde J} \to [0,\infty) ~~\text{and}~~B: 2^{\tilde J} \to [0,\infty),
$
subject to $A(\emptyset) =1$ and $B(\emptyset) =1$ and for  $ D,I,H \subseteq \tilde J$
(we set $0/0=0$ and $1/0=\infty$) subject to
$$
\frac{A(I)}{A({D})}<\infty ~~~\forall\     D  \subset I\subseteq \tilde J  \ \ \  \text{and}~~~~
\frac{B({D})}{B(H)}<\infty~~~ \forall  \  {H}\subset D\subseteq \tilde J\,. \ \ \ \  \  ~
$$
With these functions  define breakdown and repair rates as follows:
$$
q^{\Y}({D}, I)=\frac{A(I)}{A({D})},\ \ \ D\subset I\subseteq\tilde J\,,
$$
for breakdowns of nodes in non-empty set $I\setminus D$ if nodes in  $D$ are already down, and
$$
q^{\Y}(D, H)=\frac{B({D})}{B(H)},\ \ \ H\subset D \subseteq\tilde J,\,,
$$
for finishing repair of nodes in non-empty set $D\setminus H$ if nodes in $D$ are under repair.
For all other pairs $G,H\subseteq \tilde J, G\neq H,$ we set $q^{\Y}({G}, H)=0$, and
for all $D\subseteq \tilde J$ we set $q^{\Y}(D, D)= - \sum_{H\subseteq \tilde J, H\neq D} q^{\Y}(D, H)$.
\end{defn}
The generator $Q^{\Y}=(q^{\Y}(K,L): K, L \subseteq \tilde J)$ of $\Y$
is defined   for real functions
$f : 2^{\tilde J}\to \R,$ by
\begin{equation}\label{eq:bdr-generator}
    (Q^{\Y}f)(D)= \sum_{H\subset D} \frac{B(D)}{B(H)} (f(H) - f(D))
+ \sum_{I\supset D} \frac{A(I)}{A(D)} (f(I) - f(D))
\end{equation}
By evaluation of the standard local balance equations it can be seen that
$\Y$ is reversible with respect to the probability measure (with normalization constant $\hat C^{-1}$)
\begin{equation}\label{eq:pi-hat}
\hat\pi:=\left(\hat\pi( D):= \hat C^{-1} \frac{A({D})}{B({D})}, ~~ D \in 2^{\tilde J}\right) .
\end{equation}

\subsection{Rerouting}\label{sect:rerouting}
The  network process and the breakdown-repair process (availability process) interact and
we have to fix rules for the interaction regime. The general rule is:\\
{\bf (1)} Whenever a  station is broken down and under repair, service is interrupted and
the customers present there are frozen, while new customers are not admitted to this station.\\
{\bf (2)} Therefore we have to define a new routing mechanism. Examples of how to do this to obtain
explicit steady states can be found in \cite{sauer;daduna:03}[Sections 5, 6].
We describe  an  abstract ''rerouting scheme'', which encompasses the three schemes
described there.
\begin{assumption}{\sc Rerouting schemes in open  networks.}\label{defn:reroutingschemes}
Consider a Jackson network with routing matrix \eqref{eq:routing-rejection1} (with  $(r_{00}\geq 0)$)
and traffic equations for the admitted customers \eqref{eq:trafficequation-adm} where $\lambda_j=\lambda r_{0j}$.\\
When nodes in  ${D}\subseteq \tilde J$  are down,  routing is
restricted to   $ \tilde{J}_0 \setminus {D}$  and  determined by some
routing matrix
\beq\label{eq:routing-rejectionD}
R^D=[r^D_{ij}]_{i,j\in \tilde J_0\setminus D}\,.
\eeq
The associated traffic equations for the admitted customers
similar to \eqref{eq:trafficequation-adm} are
\begin{equation} \label{eq:routing-rejectionD2}
\eta^D_j=\lambda^D_j + \sum_{i\in \tilde J\setminus D} \eta^D_i r^D_{ij},\quad\quad
j\in \tilde J\setminus D,\qquad\text{with}~~\lambda^D_j := \lambda\cdot r^D_{0,j}, \quad j \in \tilde J\setminus D\,,
\end{equation}
and {\bf are assumed} to be  solved by
\begin{equation} \label{eq:routing-rejectionD3}
\eta^D_j=\eta_j,\quad j\in \tilde J\setminus D,
\end{equation}
where the $\eta_j$ are from the solution of \eqref{eq:trafficequation-adm}.
We set $R^\emptyset := R$ and, if necessary,  $\eta_j^\emptyset := \eta_j$.\\
For the  rerouting scheme with nodes in set $D$ broken down
$\lambda^D_0 := \lambda\cdot r^D_{0,0}$ is the new rejection rate.
\end{assumption}

\blemma\label{lem:reroutingschemesR}
If $(\eta^D_j, j\in \tilde J\setminus D)$ solves the traffic equations \eqref{eq:routing-rejectionD2}
 for the {admitted} customers, when nodes in $D$ are broken down and rerouting is according to
 Assumption \ref{defn:reroutingschemes},
 then with $\eta^D_{0}:= \lambda$ the vector $(\eta^D_j, j\in \tilde J_0\setminus D)$
 solves the equation $x=x\cdot R^D$.
\elemma
\bproof
\eqref{eq:routing-rejectionD2} can be written as
$\eta^D_j=\lambda\cdot r^D_{0j} + \sum_{i\in \tilde J\setminus D} \eta^D_i r^D_{ij}, \quad j \in \tilde J\setminus D\,.$
Summing over $j \in \tilde J\setminus D$ yields
$\lambda(1-r^D_{00}) = \sum_{i\in \tilde J\setminus D} \eta^D_i r^D_{i0}\,,$
which is the missing equation of $x=x\cdot R^D$, with the required solution inserted.
\eproof
\bremark\label{rem:rerouting-irred}\rm
(i) For technical simplifications we introduce the normalized probability vector $\xi^D$
associated with values
$\eta_j^D$: $\xi_j^D =\eta_j/(\sum_{k=0}^J \eta_k), j\in \tilde J_0\setminus D.$
Similarly we define $\xi=\xi^\emptyset$ associated with the $\eta_j=\eta_j^\emptyset$.\\
(ii) Prescribing rerouting by  \eqref{eq:routing-rejectionD2} is not constructive, but is not necessary for our main applications.
A detailed description of rerouting schemes which fulfill  Assumption
\ref{defn:reroutingschemes} is given in \cite{sauer:06}[Section 2].\\
(iii) When considering rerouting schemes which are used in the literature it may happen that the
rerouting chain on certain subsets $\tilde J \setminus D$ is not irreducible, for details see \cite{sauer:06}[Proof of Theorem 1.2.29].
This makes the computations more involved, but leads to the same results as those we shall present below.
\eremark

\begin{ex}\label{ex:rerouting}
The most common rerouting schemes found in the literature
which lead to explicit
stationary distributions of the network processes are (for more details see
\cite{sauer;daduna:03}[Sections 5, 6] and \cite{daduna;szekli:13}[Section 2.3])\\
(i) {\sc rs--rd with reversible routing.} This applies only for reversible  $R$.
When a customer tries to visit a  down node he is rejected, stays at his node
for another service and tries again with newly selected destination.\\
(ii) {\sc Stalling.} Whenever any node breaks down all services and arrivals are interrupted
and resumed only when all nodes are up again.\\
(iii) {\sc Skipping.} Whenever a customer wants to visit a broken down node he is
not allowed to settle
down there and  has to jump forward according to $R$ until he reaches an up-node or leaves the network.
\end{ex}

\subsection{Networks with breakdown and repair: Product formula}
\label{sect:networksJ-unrel}
A Markovian state process for an unreliable Jackson network requires that the state space
$\mathbb{N}^J$ of the Jackson network process $\X$ is supplemented by a coordinate
 $\Y$ which indicates the set of broken down stations.
Operating on these states we define a Markov process  ${{\Z}}=(\Y,\X)$ describing the degradable network with state space $\tilde\ourE =2^{\tilde J} \times  \mathbb{N}^{J}$.
States are
$\nn =( {D},\n) =( {D}, n_1,n_2,\ldots,n_J)\in  \tilde{\ourEE},$
where the  first coordinate in $\nn $ we call the availability coordinate. The interpretation is:
${D}$ is the set of servers in \emph{down status}. At node $i\in {D}$
there are $n_i$ customers waiting for server being repaired.
We denote by $\cE := 2^{\tilde{\ourEE}}$.
\begin{defn}\label{defn:unreliable-net}
{\sc The  unreliable Jackson network process} is
the Markov process ${\Z}=(Z(t),t\geq 0)$ defined by the  infinitesimal
generator (transition intensity matrix)
$Q^\Z=(q^{\Z}(\nn,\nn'):{\nn,\nn'\in\tilde \ourEE \ })$ via
\begin{eqnarray}\label{eq:generatorJ-unrel}
(Q^\Z f)({D}, n_1,n_2,\ldots,n_J)&=&
 \sum_{j \in \tilde J\setminus D} \lambda r_{0j}^D (f(D,\n + e_j) - f(\nn))
+\sum_{j \in \tilde J\setminus D} \mu_j(n_j) r_{j0}^D (f(D,\n - e_j)- f(\nn))\nonumber\\
&&+ \sum_{j \in \tilde J\setminus D} \mu_j(n_j)
\sum_{{i \in \tilde J\setminus D}}
r_{ji}^D (f(D,\n - e_j+ e_i) - f(\nn))\\
&&+ \sum_{ {H}\subset D} { B( {D})\over B({H})} (f({H},\n)-f(\nn))
+ \sum_{ {I}\supset D} { A( {I})\over A( {D})} (f({I},\n)-f(\nn)) .\nonumber
 \end{eqnarray}
\end{defn}
\begin{theorem}{\em{\sc Product form  for  Jackson networks with  breakdown and repair.}}\label{thm:pf-unreliable}
\cite{sauer;daduna:03},\\\cite{sauer:06}[Theorem 2.4.1]
Under the Assumption \ref{defn:reroutingschemes}, if ${\Z}=(\Y,\X)$ is ergodic then the steady state is with $\pi$ from \eqref{eq:steadystateJ} and
$\hat \pi$ from \eqref{eq:pi-hat}  of the product form:
~For $({D}, \: n_1,\dots,n_J) \in \tilde\ourEE$
\beq\label{eq:pf-unreliable}
\tilde\pi(D, \: n_1, n_2,\dots, n_J)= \hat \pi(D) \cdot  \pi(n_1, n_2,\dots, n_J)
= \hat C^{-1} \frac{A(D)}{B(D)} \cdot
\prod_{j=1}^J \left(C_j^{-1}\prod_{i=1}^{n_j} \frac{\eta_j}{\mu_j(i)}\right)\,.
\eeq
Here $\eta = (\eta_1, \dots, \eta_J)$ is the solution of the traffic equation
 \eqref{eq:trafficequation-adm} for {admitted} customers when all nodes are up,
$C_j$ is the normalization constant for the local queue length process at node $j$. We denote
$C =
 \prod_{j=1}^J C_j\,.$
\end{theorem}
\eqref{eq:pf-unreliable}  is proved in \cite{sauer;daduna:03}  for more general breakdown and repair schemes:
Breakdown and repair rates may depend on the load (queue lengths) of nodes.
The question whether in this framework results similar to those in the following sections can be
derived is still open and part of our ongoing research.\\
Theorem \ref{thm:pf-unreliable} is not covered by the results for networks in a random environment in
\cite{zhu:94} and \cite{economou:05}. In both papers  it assumed that under different environment
states the ratio ''local arrival rate/local service rate'' is independent of the environment status.
 This is obviously not the case in our systems.

\section{One step  correlation}
\label{sect:1stepcorr}
 Recall $\lambda_j = \lambda\cdot r_{0j}, j=1,\dots,J,$
and  $\lambda_1+\ldots +\lambda_J = \lambda(1-r_{00})$ and that we therefore consider
only admitted customers even if all nodes are up. We will not mention this further.
For the network process $\Z$ with generator $Q^{\Z}$
and stationary distribution $\tilde\pi$ consider \emph{one step}
correlation expressions
\begin{equation}\label{eq:corr1}
\langle f,Q^{\Z} g\rangle_{\tilde\pi}. 
\end{equation}
If  $f=g$, then  \eqref{eq:corr1} is (the negative of) a quadratic form, because
$-Q^{\Z}$ is positive definite. \eqref{eq:corr1} occurs in the definition of
Cheeger's constant because division of \eqref{eq:corr1} with $f=g$ by $\langle f,f\rangle_{\tilde\pi}$ yields Rayleigh quotients. It also occurs in the definition of the corresponding Dirichlet form.
This is helpful to bound the second largest eigenvalue of $Q^{\Z}$ and to prove the Poincare inequality for the corresponding  Markovian process, see e.g., \cite{chen:04}.
Furthermore, \eqref{eq:corr1} can be utilized to determine the asymptotic variance of costs or performance measures associated with Markovian network processes  and to compare the  asymptotic variances of two such related processes. It is possible to compare the correlations for $\Z$
with that of the related process ${\Z'}$ with the same  stationary distribution $\tilde\pi$, using
$\langle f,Q^{\Z} g\rangle_{\tilde\pi}-\langle f,Q^{\Z'} g\rangle_{\tilde\pi},$
which will be given explicitly in Section \ref{sect:comparison-corr}.
Because we are dealing with processes having bounded generators, properties connected with \eqref{eq:corr1} can be turned into properties of
$\langle f,Id +\varepsilon Q^{\Z} g\rangle_{\tilde\pi} = E_{\tilde\pi} (f(Z_0)g(Z_\tau))$
where $\varepsilon>0$ is sufficiently small such that
$Id +\varepsilon Q^{\Z}$ is a stochastic matrix and $\tau$ is exponentially distributed.
This enables one to directly apply discrete time methods to characterize
properties of continuous time processes.
\subsection{ Correlation formulas}
\label{sect:corrformula}
Due to the product form steady state distribution of $\Z$ the one step correlation
$\langle f,Q^{{\Z}} g\rangle_{\tilde\pi}$ splits immediately into two terms having an intuitive interpretation: The sum of weighted
one step conditional correlations\\  $(i)$ of the environment process $Y$, and
 $(ii)$ of network processes, which for a fixed time point
 seem to behave conditionally independent of the environment.\\
As will be seen, it is illuminating to define for all $D\subseteq \tilde J$
the generators $Q^{\X_{\tilde J\setminus D}}$ of ''synthetic subnetworks'' on node set
$\tilde J \setminus D$ with overall arrival rate $\lambda$, service rates from
Definition \ref{defn:unreliable-net}, and routing matrix $R^D$.\\
The proofs of Propositions \ref{prop:corrformula-split} and \ref{prop:corr-formula}
are postponed to Section \ref{sect:proofs}.
\bprop\label{prop:corrformula-split}{\sc Splitting formula.}
 For unreliable Jackson network processes ${\Z}=(\Y,\X)$
the one step correlations splits as follows
\beqaro
&&\langle f,Q^{{\Z}} g\rangle_{\tilde\pi}
= \sum_{\n\in \ourN^J} \pi(\n)
\left\{\sum_{D\subseteq {\tilde J}} \hat\pi(D)
f(D,\n) \left( Q^{\Y} g(\cdot,\n)\right)(D)\right\}
+ \sum_{D\subseteq {\tilde J}} \hat\pi(D)
\sum_{\n_D\in \ourN^{|D|}} \pi_D(\n_D)\nonumber\\
&&\qquad\qquad\qquad\qquad\qquad\cdot \left\{
\sum_{\n_{\tilde J\setminus D}\in \ourN^{J-|D|}}\pi_{\tilde J\setminus D}(\n_{\tilde J \setminus D})
f(D,(\n_D,\n_{\tilde J \setminus D}))
\left(
Q^{\X_{\tilde J\setminus D}}g(D,(\n_D,{(\cdot)}_{\tilde J \setminus D}))
\right)(\n_{\tilde J \setminus D})
\right\}
\eeqaro
\eprop
The  next correlation formula  will yield remarkable simplifications when used for differences.
\begin{prop}\label{prop:corr-formula}
{\sc One-step correlation  formula.}
 For unreliable Jackson network processes ${\Z}=(\Y,\X)$
with $\xi^D$ the probability solution of
$x^D=x^D\cdot R^D$ (when nodes in $D$ are down) holds
\begin{eqnarray}
  \langle f,Q^{{\Z}} g\rangle_{\tilde\pi}
  = \sum_{D\subseteq {\{1,\dots,J\}}} \sum_{\n\in \ourN^J} \tilde\pi(D,\n)\ \   
 \left\{\sum_{H\subset D} q^{\Y}(D,H) f(D,\n) g(H,\n)\right.
  \left.+
 \sum_{I\supset D} q^{\Y}(D, I) f(D,\n) g(I,\n)
\right\}\nonumber\\
+ \sum_{D\subseteq {\{1,\dots,J\}}} \sum_{\n\in \ourN^J}\tilde\pi(D,\n) \frac{\lambda}{\xi^D_0}
\left\{\sum_{j\in \tilde J_0 \setminus D} \sum_{i\in \tilde J_0 \setminus D}
\xi^D_j r^D_{ji} f(D,\n+e_j) g(D,\n+e_i)\right\} \qquad\quad\qquad\label{corrnet4}\\
- \sum_{D\subseteq {\{1,\dots,J\}}} \sum_{\n\in \ourN^J} \tilde\pi(D,\n) f(D,\n) g(D,\n)
\left\{\sum_{H\subset D} q^{\Y}(D, H)  + \sum_{I\supset D} q^{\Y}(D, I) 
+\lambda + \sum_{j\in \tilde J \setminus D}
\mu_j(n_j) \right\}.\nonumber
\end{eqnarray}
\end{prop}

\subsection{Comparison of one step correlations}\label{sect:comparison-corr}

The following formulas for differences of one step correlations will give additional insight into  various properties of networks, for example to speed of convergence or asymptotic variance. They display how, e.g., the routing and
the breakdown and repair affects correlations in networks.
\begin{theorem}\label{thm:corrgeneral}{\sc Changes of routing
behaviour of the customers.}
Suppose ${\Z}=(\Y,\X)$ is an ergodic unreliable Jackson network process  with a
routing matrix $R$ and $\Z'=(\Y',\X')$ is another Jackson network
process having the same arrival and service intensities and failure-repair rates but with routing
matrix $R'=[ r'_{ij}]_{i,j=0,1,\dots,J}$, such that the
solutions of the traffic equation derived from $R$ and for $ R'$
coincide (denoted by $\eta$).  Assume that  both networks  follow some rerouting mechanism for which  the Assumption \ref{defn:reroutingschemes} holds.
Then for arbitrary real functions $f,g\in L^2$
$$\langle f,Q^{{\Z}} g\rangle_{\tilde\pi}-\langle f,Q^{\Z'} g\rangle_{\tilde\pi}=
 E_{\tilde\pi}\left[ \frac{\lambda}{\xi^{Y_t}_0}\left(tr\big(W^{g,f}(Y_t,X_t)\cdot diag(\xi^{Y_t})\cdot(R^{Y_t}-
{R'}^{Y_t})\big)\right)\right], $$
where $\xi^D$ is the probability solution of $x^D = x^D \cdot R^D$,
$tr(A)$ denotes trace of $A$,  and 
$$W^{g,f}(D,\n)=[g(D,\n+e_i)f(D,\n+e_j)]_{i,j\in {\tilde J}_0\setminus D}\,.
$$
\end{theorem}
\bproof Because the external arrival streams are the same, and the traffic equations have the same solution $\eta$,
and the rerouting mechanisms have property \eqref{eq:routing-rejectionD3}, for any availability status $D$ the
 rerouting schemes on $\tilde J\setminus D$ have the same solution of the traffic equation.
It follows from Lemma \ref{lem:reroutingschemesR} that for all $D$ the  probability solution of the
equations $x^D = x^D \cdot R^D$ and $x^D = x^D \cdot R'^D$ are in both systems the same.
Because of $q^{\Y} = q^{\Y'}$ we  have from Proposition
\ref{prop:corr-formula} the reduction
\beqaro
&&\langle f,Q^{{\Z}} g\rangle_{\tilde\pi}-\langle f,Q^{\Z'} g\rangle_{\tilde\pi}= \sum_{D\subseteq {\{1,\dots,J\}}} C^{-1}\hat \pi(D) \sum_{\n_D\in \ourN^{|D|}}
\prod_{\ell\in D} \prod_{i=1}^{n_\ell} \left(\frac{\eta_\ell}{\mu_\ell(i)}\right)\\
&&\quad\quad\qquad\qquad\qquad\cdot \sum_{\n_{\tilde J\setminus D}\in \ourN^{J-|D|}}
\prod_{\ell\in \tilde J\setminus D} \prod_{i=1}^{n_\ell} \left(\frac{\eta_\ell}{\mu_\ell(i)}\right)
\left[
\frac{\lambda}{\xi^D_0}
\left\{\sum_{j\in \tilde J_0 \setminus D} \sum_{i\in \tilde J_0 \setminus D}
\xi^D_j r^D_{ji} f(D,\n+e_j) g(D,\n+e_i)\right\}\right.\\
&&\quad\qquad\qquad\qquad\qquad\qquad\qquad\qquad\qquad\qquad\qquad\left. - \frac{\lambda}{\xi^D_0}
\left\{\sum_{j\in \tilde J_0 \setminus D} \sum_{i\in \tilde J_0 \setminus D}
\xi^D_j {r'}^D_{ji} f(D,\n+e_j) g(D,\n+e_i)\right\}
\right]
\eeqaro
We interpret in the last two lines for fixed $D$ and $\n_D$ and $i,j\in {\tilde J_0 \setminus D}$ the expressions
\beqo
f(D,\n+e_j)=: f(D,(\n_D, \n_{\tilde J \setminus D}+e_j))\quad\text{and}\quad
 g(D,\n+e_i) =: g(D,(\n_D, \n_{\tilde J \setminus D}+e_i))
\eeqo
as functions of $\n_{\tilde J \setminus D}$ only, and see that the resulting expressions have
exactly the structure of the functions dealt with in Proposition 4.1 of \cite{daduna;szekli:08}.
After renormalization of the densities
$\prod_{\ell\in \tilde J\setminus D} \prod_{i=1}^{n_\ell} \left(\frac{\eta_\ell}{\mu_\ell(i)}\right)$,
which  in fact results in conditioning on $\{Y(t)=D,X_{D}(t)=\n_D\}$, we obtain
$\langle f,Q^{{\Z}} g\rangle_{\tilde\pi}-\langle f,Q^{\Z'} g\rangle_{\tilde\pi}$
\beqaro
&=& \sum_{D\subseteq {\{1,\dots,J\}}}  \sum_{\n_D\in \ourN^{|D|}} P(\{Y(t)=D,X_{D}(t)=\n_D\})\\
&&\qquad\qquad\cdot E_{\tilde \pi}\left[\frac{\lambda}{\xi^D_0} \left(tr(W^{g,f}(D,(\n_D,X_{\tilde J_0 \setminus D}))
\cdot diag\xi^D\cdot(R^D-{R'}^D))\right)|\{Y(t)=D,X_{D}(t)=\n_D\}\right]\,,
\eeqaro
and deconditioning eventually finishes the proof.
\eproof

\begin{theorem}\label{thm:corrbdr}{\sc Changes of breakdown and repair mechanisms.}
Suppose ${\Z}=(\Y,\X)$ is an ergodic unreliable Jackson network process  with a
routing matrix $R=[ r_{ij}]_{i,j=0,1,\dots,J}$ and $\Z'=(\Y',\X')$ is another Jackson network
process having the same arrival and service intensities, and with the same routing regime, described
by  $R$ and rerouting fulfilling Assumption \ref{defn:reroutingschemes}.\\
The breakdown-repair process  for $\Z$ is given in Definition \ref{defn:queueindependentbreakdown} and
is for $\Z'$  defined similarly via
$A',~B': 2^{\tilde J} \to [0,\infty),$
subject to the restrictions  indicated there.
Then the  breakdown and repair rates for $\Y'$ are:
$$
q^{\Y'}({D}, I)=\frac{A'(I)}{A'({D})},\ \ \ D\subset I\subseteq\tilde J,\qquad
\text{and}~~q^{\Y'}(D, H)=\frac{B'({D})}{B'(H)},\ \ \ H\subset D \subseteq\tilde J,\,.
$$
The processes $\Y$ and $\Y'$ are Markov  with generators $Q^{\Y}=(q^{\Y}(K,L): K, L \subseteq \tilde J)$ of $\Y$
and $Q^{\Y'}=(q^{\Y'}(K,L): K, L \subseteq \tilde J)$ of $\Y'$ as defined in \eqref{eq:bdr-generator} for
$\Y$ and similar for $\Y'$.

Assume that the stationary distributions of $\Y$ and $\Y'$ are identical, denoted by
\begin{equation*}
\hat\pi:=\left(\hat\pi( D):= \hat C^{-1}\frac{A({D})}{B({D})} = \hat C'^{-1}\frac{A'({D})}{B'({D})}, ~~ D \in 2^{\tilde J}\right)\,.
\end{equation*}
  Then for arbitrary real functions $f,g: \tilde\ourEE\to \R$ holds
  \begin{equation*}\label{eq:bdr-differences1}
\langle f,Q^{{\Z}} g\rangle_{\tilde\pi}-\langle f,Q^{\Z'} g\rangle_{\tilde\pi}
= E_{\pi}
\left[ \langle f(\circ,X_t),Q^{{\Y}} g(\cdot,X_t)(\circ)\rangle_{\hat\pi}-\langle f(\circ,X_t),Q^{\Y'} g(\cdot,X_t) (\circ)\rangle_{\hat\pi}
\right]=
  \end{equation*}
$$
=E_{\pi}
\left[ \langle f(\circ,X_t),(Q^{{\Y}}-Q^{\Y'}) g(\cdot,X_t)(\circ)\rangle_{\hat\pi}
\right].
$$
\end{theorem}
\bproof Interchanging summations, regrouping terms, and exploiting the product form structure
of the state distributions (which are identical for
${\Z}=(\Y,\X)$ and $\Z'=(\Y',\X')$) in \eqref{corrnet4}
we obtain
\beqaro
\langle f,Q^{{\Z}} g\rangle_{\tilde\pi}
&=& \sum_{\n\in \ourN^J} \pi(\n) \sum_{D\subseteq {\{1,\dots,J\}}} \hat\pi(D)\ \   
 \left\{\sum_{H\subset D} q^{\Y}(D,H) f(D,\n) g(H,\n)
 + \sum_{I\supset D} q^{\Y}(D, I) f(D,\n) g(I,\n)
\right\}\nonumber\\
&&- \sum_{\n\in \ourN^J} \pi(\n) \sum_{D\subseteq {\{1,\dots,J\}}}  \hat\pi(D) f(D,\n) g(D,\n)
\left\{\sum_{H\subset D} q^{\Y}(D, H)  + \sum_{I\supset D} q^{\Y}(D, I)\right\}\\
&&+ \sum_{D\subseteq {\{1,\dots,J\}}} \sum_{\n\in \ourN^J}\tilde\pi(D,\n) \frac{\lambda}{\xi^D_0}
\left\{\sum_{j\in \tilde J_0 \setminus D} \sum_{i\in \tilde J_0 \setminus D}
\xi^D_j r^D_{ji} f(D,\n+e_j) g(D,\n+e_i)\right\}\\
&&- \sum_{D\subseteq {\{1,\dots,J\}}} \sum_{\n\in \ourN^J} \tilde\pi(D,\n) f(D,\n) g(D,\n)
\left\{\lambda + \sum_{j\in \tilde J \setminus D}\mu_j(n_j)\right\}.
\eeqaro
For $\langle f,Q^{\Z'} g\rangle_{\tilde\pi}$ the last two lines  in
the respective formula are identical to those in the displayed formula.
The difference therefore is  $\langle f,Q^{{\Z}} g\rangle_{\tilde\pi}-\langle f,Q^{\Z'} g\rangle_{\tilde\pi}=$
\beqaro
&=& \sum_{\n\in \ourN^J} \pi(\n) \sum_{D\subseteq {\{1,\dots,J\}}} \hat\pi(D)\ \  
 \left\{\sum_{H\subset D} q^{\Y}(D,H) f(D,\n) g(H,\n)
  \right.+ \sum_{I\supset D} q^{\Y}(D, I) f(D,\n) g(I,\n)\\
&&  - \sum_{\n\in \ourN^J} \pi(\n) \sum_{D\subseteq {\{1,\dots,J\}}}  \hat\pi(D) f(D,\n) g(D,\n)
\left\{\sum_{H\subset D} q^{\Y}(D, H)  + \sum_{I\supset D} q^{\Y}(D, I)\right\}\\
&&- \sum_{\n\in \ourN^J} \pi(\n) \sum_{D\subseteq {\{1,\dots,J\}}} \hat\pi(D)\ \   
 \left\{\sum_{H\subset D} q^{\Y'}(D,H) f(D,\n) g(H,\n)
 + \sum_{I\supset D} q^{\Y'}(D, I) f(D,\n) g(I,\n)
\right\}\nonumber\\
&&+ \sum_{\n\in \ourN^J} \pi(\n) \sum_{D\subseteq {\{1,\dots,J\}}}  \hat\pi(D) f(D,\n) g(D,\n)
\left\{\sum_{H\subset D} q^{\Y'}(D, H)  + \sum_{I\supset D} q^{\Y'}(D, I)\right\}
\eeqaro
For fixed   $\n$ we interpret  $f(D,\n)$ and $g(D,\n)$ as functions of $D$ parametrized by  $X_t(\omega) = \n$.
This leads to
\beqar\label{eq:corrbdr-1}
&&\langle f,Q^{{\Z}} g\rangle_{\tilde\pi}-\langle f,Q^{\Z'} g\rangle_{\tilde\pi}
= \sum_{\n\in \ourN^J} \pi(\n)
\left[\langle f(\circ,\n),Q^{{\Y}} g(\cdot,\n)(\circ)\rangle_{\hat\pi}-\langle f(\circ,\n),Q^{\Y'} g(\cdot,\n)(\circ)\rangle_{\hat\pi}
\right]\nonumber\\
&&= E_{\pi}
\left[ \langle f(\circ,X_t),Q^{{\Y}} g(\cdot,X_t)(\circ)\rangle_{\hat\pi}-\langle f(\circ,X_t),Q^{\Y'} g(\cdot,X_t) (\circ)\rangle_{\hat\pi}
\right].
\eeqar
\eproof
Theorems \ref{thm:corrgeneral} and \ref{thm:corrbdr} state a reduction of complexity:
Comparing operators and resulting one step correlations via functions on infinite  space
is reduced to comparing
matrix operators via functions on finite  space.
The theorems  are valid for  all square integrable
functions $f,g$ on $\tilde \ourE$. This opens the way to compare  correlations for
multidimensional marginals in time of the network processes with unreliable nodes according to
concordance ordering utilizing the abstract setting for Markov processes
of \cite{daduna;szekli:08}[Theorem 5.2].\\
For a  concise notation we introduce the standard difference operators for functions on $\ourN^J$.
For all $f\in L^2$ and all $j=0,1,\dots,J,$ we define (recall $e_0$ is the zero vector)
\begin{equation*}
    {\cal D}_j f: \tilde\ourEE\to \R,\quad (D,\n) \to {\cal D}_j f(D,\n) := f(D,\n + e_j) - f(D,\n),
\end{equation*}
and
$$
 {\cal D} f: \tilde\ourEE\to \R^{J+1-|D|}, \quad  (D,\n) \to ({\cal D}_j f(D,\n), j\in \tilde J_0\setminus D).
$$
That way we can treat ${\cal D}f(D,\n)$ as a vector of the dimension corresponding to the size of $D$, and the corresponding routing matrices $R^D$ as operators on it. Moreover, it is possible to
consider the corresponding scalar products generated by invariant vectors $\xi^D$, and write the formula for the difference of one step correlations (which are scalar products with respect to invariant measures for the network process ${\tilde\pi}$) in terms of scalar products with respect to invariant measures $\xi^D$ for the routing processes.
\begin{cor}
For unreliable Jackson network processes ${\Z}=(\Y,\X)$
with $\xi^D$ the probability solution of the
 equation $x^D=x^D\cdot R^D$ and $\underline\mu _D=\sum_{j\in \tilde J \setminus D}
\mu_j(n_j)$ the total service rate for nodes which are up,
we have
\begin{eqnarray*}\label{correlation3}
  \langle f,Q^{{\Z}} g\rangle_{\tilde\pi}&
= &E_{\tilde\pi} \left[\frac{\lambda}{\xi^{Y_t}_0}
\left\langle ({\cal D}+Id)f(Y_t,X_t),R^{Y_t}
({\cal D}+Id)g(Y_t,X_t)\right\rangle_{\xi^{Y_t}}\right]\\
&+&E_{\pi}
\left[ \langle f(\circ,X_t),Q^{{\Y}} g(\cdot,X_t)(\circ)\rangle_{\hat\pi}
\right]
-E_{\tilde\pi}\left[f(Y_t,X_t)g(Y_t,X_t)(\lambda+\underline\mu_{ Y_t})\right].
\end{eqnarray*}
\end{cor}
\bproof
Take  Proposition \ref{prop:corr-formula} \eqref{corrnet4} , condition as in
Theorem \ref{thm:corrgeneral}, and insert the suitable difference operators.
\eproof
We can now reformulate the result of  Theorem \ref{thm:corrgeneral} in a more compact
form which immediately relates our results  to methods dealt with in optimizing MCMC simulation.
\begin{cor}\label{cor:corrgeneral*}
 For unreliable Jackson network processes ${\Z},\Z'$ as in
Theorem  \ref{thm:corrgeneral} we have
\begin{equation*}
\langle f,Q^{{\Z}} g\rangle_{\tilde\pi}-\langle f,Q^{\Z'} g\rangle_{\tilde\pi}
= E_{\tilde\pi} \left[\frac{\lambda}{\xi^{Y_t}_0}
\left\langle ({\cal D}+Id)f(Y_t,X_t),\big(R^{Y_t}-{{R'}^{Y_t}}\big)
({\cal D}+Id)g(Y_t,X_t)\right\rangle_{\xi^{Y_t}}\right]\,.
\end{equation*}
\end{cor}
\section{Applications}
\label{sect:applications}
The proofs of the applications below follow the ideas
which are used to prove the theorems in Section \ref{sect:comparison-corr}. We
therefore give only hints to perform the proofs here.
\subsection{Comparison of spectral gaps}\label{sect:comparison-gaps}
Let ${\Z}$ be a continuous time homogeneous ergodic Markov process with
stationary probability $\tilde\pi$
and generator $Q^{\Z}$. Let $\tilde\pi (f)=\int_{\tilde\ourE} f(x) \tilde\pi(dx)$
The spectral gap of ${\Z}$, resp. $Q^{\Z}$ is
\begin{equation*}
Gap(Q^{\Z}) = \inf\{\langle f,-Q^{\Z} f\rangle_{\tilde\pi}:f\in L^2(\tilde\ourE,\tilde\pi), \tilde\pi(f)=0,
~\langle f,f\rangle_{\tilde\pi}=1\}.
\end{equation*}
The spectral gap determines for ${\Z}$ the speed of convergence to equilibrium
$\tilde\pi$ in $L^2(\tilde\ourE,\tilde\pi)$ with norm $\| \cdot\|_{\tilde\pi}$: \,$Gap(Q^{\Z})$
is the largest number $\Delta$ such that for the transition semigroup
$P=(P_t:t\geq 0)$ of ${\Z}$ holds
\begin{equation*}
\| P_tf - \tilde\pi(f)\|_{\tilde\pi} \leq e^{-\Delta t}\| f - \tilde\pi(f)\|_{\tilde\pi}
~~~\forall f\in L^2(\tilde\ourE ,\tilde\pi).
\end{equation*}
It should be noted that one has to be careful which class of functions is used for the definition of
spectral gap. For a discussion and more references see the introduction of
\cite{lorek;szekli:14}.\\
We utilize the following orderings to compare routings, failure processes and then correlations, see \cite{peskun:73}.
\begin{defn}\label{defn:peskun}
Let $R=[  r_{ij}]$ and $ R'=[r'_{ij}]$ be transition
matrices on a finite set $\mathbb E$ such that $\xi R=\xi  R'=\xi$ for a probability vector $\xi$.\\
$R'$  is smaller than $R$ in the positive semidefinite order,
$R ' \prec_{pd}  R$, if  the matrix $R- R' $ is positive semidefinite.\\
$ R'$  is smaller than $R$ in the Peskun order,  $ R' \prec_P R$,
if for all $j, i\in \mathbb E$ with $i\neq j$  holds $ r'_{ji}\leq  r_{ji}$.
\end{defn}
Peskun used the latter order to compare reversible transition matrices with the
same stationary distribution  and their asymptotic variance, and Tierney
\cite{tierney:98} has shown (in a more general setting, i.e. using operators rather than matrices) that the main property used in the proof of Peskun,
namely  that ''$R \prec_P R'$ implies $ R' \prec_{pd} R$'',
holds without reversibility assumptions.
\bex\label{ex:peskun1}\rm
For any transition matrix $R=[  r_{ij}]$ and $R'=Id$ (of the same dimension as $R$) holds
$R' \prec_P R,$ so the family of transition matrices of  fixed dimension has a (unique) minimal
element. Therefore {\sc stalling}
from Example \ref{ex:rerouting} is an extremal (re-)routing scheme, because it utilizes $Id$
in case of any breakdown.
\eex
\begin{prop}\label{prop:pdorder-nets}
 Consider two ergodic unreliable Jackson networks  with state processes  ${\Z}$ and ${\Z}'$
 and with the same arrival and service
intensities, and the same failure-repair rates.
Assume that  the  equations
$x=x\cdot R$ and $x=x\cdot R'$ have the same normalized solution $\xi$, and the Assumption
 \ref{defn:reroutingschemes} holds, i.e. both networks follow some  rerouting mechanism according to \eqref{eq:routing-rejectionD2}
with the property \eqref{eq:routing-rejectionD3}.\\
If   $R^D\prec_{pd}R'^D$  for all $D$ (also for $D=\emptyset$), then
\begin{equation*}
Gap(Q^{{\Z}')}\le Gap(Q^{\Z}).
\end{equation*}
\end{prop}
\bproof
From Corollary  \ref{cor:corrgeneral*} we have for all $f\in L^2$
\begin{equation*}\label{eq:independent-j}
\langle f,-Q^{{\Z}} f\rangle_{\tilde\pi}-\langle f,-Q^{\Z'} f\rangle_{\tilde\pi}
= E_{\tilde\pi} \left[\frac{\lambda}{\xi^{Y_t}_0}
\left\langle ({\cal D}+Id)f(Y_t,X_t),\big(R'^{Y_t}-{{R}^{Y_t}}\big)
({\cal D}+Id)f(Y_t,X_t)\right\rangle_{\xi^{Y_t}}\right],
\end{equation*}
and from the product formula we rewrite this formula as
\begin{equation*}
\langle f,-Q^{{\Z}} f\rangle_{\tilde\pi}-\langle f,-Q^{\Z'} f\rangle_{\tilde\pi}
= E_{\pi}E_{\hat\pi} \left[\frac{\lambda}{\xi^{Y_t}_0}
\left\langle ({\cal D}+Id)f(Y_t,X_t),\big(R'^{Y_t}-{{R}^{Y_t}}\big)
({\cal D}+Id)f(Y_t,X_t)\right\rangle_{\xi^{Y_t}}\right]=
\end{equation*}
\begin{equation*}
= E_{\pi}\sum_D\hat\pi(D) \left[\frac{\lambda}{\xi^{D}_0}
\left\langle ({\cal D}+Id)f(D,X_t),\big(R'^{D}-{{R}^{D}}\big)
({\cal D}+Id)f(D,X_t)\right\rangle_{\xi^{D}}\right].
\end{equation*}
From  $\prec_{pd}$ ordering of routings for all $f\in L^2$ therefore
$
\langle f,-Q^{{\Z}} f\rangle_{\tilde\pi}\ge\langle f,-Q^{\Z'} f\rangle_{\tilde\pi}
$
holds.
\eproof
\begin{rem}\label{rem:stalling-extremal}
If $R\prec_{pd} R'$ and the rerouting for
$\Z'$ is by stalling, then by the extremal property of $Id$ for $\prec_P$
 the assumptions of Proposition \ref{prop:pdorder-nets} are fulfilled.
\end{rem}
Computing spectral gaps for Markov processes with  multidimensional state space is challenging,
in many cases nearly impossible. Exceptions are multidimensional independent birth-death processes,
because for birth-death processes explicit results are known, see e.g. \cite{doorn:02}, and Liggett
has proved that the gap of independent processes is the minimum of the spectral gap of the marginal
processes \cite{liggett:89}[Theorem 6.2].
We will show that the gap of the joint queue length network process ${\Z}$ (with unreliable nodes)
can be bounded from below by the gap of a network process consisting of identical breakdown-repair process and  related  multidimensional birth-death process with
{\bf conditionaly independent} components.

\begin{prop}\label{prop:gap-bd-abstract}
Consider an ergodic Jackson network process  ${\Z}$ with unreliable servers as in
Theorem \ref{thm:pf-unreliable}.
 Assume that for all $D\subseteq \{1,2,\dots,J\}$ with $\hat\pi(D)>0$ the  routing matrix
 $R^D=[r_{ij}^D]_{i,j\in \tilde J_0 \setminus D}$ has strict positive entrance and departure
 probabilities ($r_{0i}^D>0, r_{i0}^D>0$)  for every node $i\in \tilde J \setminus D$.\\
 Assume further that for all $D\subseteq \{1,2,\dots,J\}$ the rerouting $R^D$ fulfills
 overall balance for all network nodes which are up, i.e.
 \begin{equation}\label{eq:detbalance1}
 \eta_j^D \sum_{i\in \tilde J \setminus D}  r_{ji}^D= \sum_{i\in \tilde J \setminus D} \eta_i^D r_{ij}^D,\quad
 \forall j\in \tilde J \setminus D.
 \end{equation}
Then there exists  an ergodic Jackson network process process $\Z'$ with unreliable servers as
in Theorem \ref{thm:pf-unreliable} with the same stationary distribution $\tilde \pi$
as $\Z$, such that
$$Gap(Q^{\Z'})\le Gap(Q^{\Z})\,.$$
The nodes of $\Z'$ are perturbed by a common breakdown-repair regime identical to that of $\Z$,
and for any given set $D$ of broken down nodes the joint network process on  $\tilde J\setminus D$
consists of conditionally independent birth-death processes,
and the coordinate birth and death processes on the $i$-th coordinate have birth rate
$\lambda r^D_{0i}$ and  state dependent death rate $\mu_i(n_i)r^D_{i0}$.
\end{prop}
\bproof
In order to obtain a lower bound for the spectral gap using  birth and death processes the idea  is to allow in the comparison network for each node $i,$ which is up, that any customer who enters node $i$ from the external source after being served
only to feed back (possibly iteratively) to node $i$ or to depart from the network. This results
in updating the service rates suitably. Recall that the original network processes and the
constructed comparison processes are additionally perturbed by the same failure mechanism.\\
Consider the situation when nodes in $D$ are down.
Directly from the generator of the network process $\Z$ in \eqref{eq:generatorJ-unrel},
it is clear that after reducing movements inside the  network, and allowing only for movements into the network from outside, or from the network into outside, or feedback, we get as long as the
reliability level $D$ does not change transitions for changing the queue lengths which look
identical as those of the generator
of independent birth and death  processes such that on the $i$-th coordinate  the birth rate equals $\lambda r^D_{0i}$ and the state dependent death rate equals
$\mu_i(n_i)r^D_{i0}, i\in \tilde J\setminus D$.\\
Now, in order to be able to apply a formula for differences of one step correlations we have to show
on every reliability level $D$
that such a modification is possible within a class of networks with extended routings having the same stationary solution. For this reason we need the assumption on overall balance (\ref{eq:detbalance1}).\\
More precisely, we define $R'^D$  by $r'^D_{i0}=r^D_{i0}$, $r'^D_{0i}=r^D_{0i}$, for all $i$,  $r'^D_{ij}=0$ for $j\neq i,\ i,j\in \tilde J\setminus D$, and $r'^D_{ii}=1-r^D_{i0}$ for
$i\in \tilde J\setminus D$. With the routing $R'^D$, the network process $\Z'$  (when nodes in $D$ are down) develops as a vector of independent birth and death processes for the up nodes
in $\tilde J\setminus D$.\\
For $j\in \{1,\dots,J\}$ let $\eta'^D_j$ be the solution of the traffic equations for $R'^D$.
We have directly
\beq\label{eq:routing-birthdeath1}
{\eta}'^D_j= \lambda r^D_{0j}+ {\eta}'^D_j  r'^D_{jj},~~~~j\in \tilde J\setminus D\,,
\eeq
and the solution of this system is uniquely defined.  We show that \eqref{eq:routing-birthdeath1} is solved
(for each $D$) by $({\eta}^D_j, j\in \tilde J\setminus D)$
as well. Inserting this
into \eqref{eq:routing-birthdeath1} we obtain with $r'^D_{jj}=1-r^D_{j0}$
for $j\in \tilde J\setminus D$
\beqo
{\eta}^D_j= \lambda r^D_{0j}+ {\eta}^D_j (1-r^D_{j0})
= \lambda r^D_{0j}+ {\eta}^D_j \sum_{i\in \tilde J \setminus D} r^D_{ji}
= \lambda r^D_{0j}+  \sum_{i\in \tilde J \setminus D} {\eta}^D_i r^D_{ij}\,,
\eeqo
which is the traffic equation when nodes in $D$ are down in $\Z$ and has the
unique solution $({\eta}^D_j, j\in \tilde J\setminus D)$.
The last step is: First observe that $({\eta}^D_j, j\in \tilde J\setminus D)$
is by Assumption \ref{defn:reroutingschemes} the restriction of $({\eta}_j, j\in \tilde J)$
(the solution of the
traffic equation with all nodes in $\Z$ up), to $\tilde J \setminus D$, and, secondly,
consider the above constructed  system of independent birth-death processes for the
reliability level $\emptyset$ as the comparison system when all nodes are up with routing $R'$,
and the $R'^D$ as rerouting scheme for this network on reliability level $D.$\\
Because ${\eta_j}={\eta'_j}$ and ${\eta}^D_j = {\eta}'^D_j$ for all $j\in \tilde J \setminus D$ and all $D$, the ${\eta}'^D_j$ are obtained by
restricting  $({\eta}_j':j\in \tilde J)$ (the solution of the
traffic equation with all nodes in $\Z'$ up), to $\tilde J \setminus D$.\\
Note that for all $D$ holds
$R'^D\prec_{P}R^D,$
therefore $R^D\prec_{pd}R'^D,$ and the result follows from Proposition \ref{prop:pdorder-nets}.
\eproof
\begin{rem}\label{rem:special-birth-death}
The overall balance \eqref{eq:detbalance1} of the rerouting $R^D$ for all network nodes which are up
is satisfied whenever $R$ is a reversible transition matrix.\\
If the overall balance \eqref{eq:detbalance1} holds for $\tilde J$ only (i.e., for $D=\emptyset$),
and the rerouting for
$\Z'$ is by stalling, the extremal property of $Id$ for $\prec_P$ leads to the same conclusion
as the assumptions of Proposition \ref{prop:gap-bd-abstract}.
\end{rem}
Remark: From irreducibility of $\Z$, for  $i, j\in \tilde J, i\neq j$  we obtain
from \eqref{eq:detbalance1}
 $\eta_j  r_{ji}^{\tilde J\setminus \{i,j\}}= \eta_i r_{ij}^{\tilde J\setminus \{i,j\}},$
but this does not mean that the matrix $R$ is reversible, because  $r_{ij}^{{\tilde J\setminus \{i,j\}}}\neq r_{ij}$ may hold.

The lower bound  $Gap(\Z')$ in the previous statements is of interest, because it has
 constitutive processes with conditionally independent coordinates.
From \cite{liggett:89}[Theorem 2.6] it is known, that the gap of a process with independent
coordinates is the minimum of the gaps of the coordinate processes. Unfortunately enough,
this theorem does not apply here
directly, because the coordinate birth-death processes are controlled by the
common breakdown-repair process. However, the comparison result of Proposition \ref{prop:pdorder-nets} can be used to obtain upper bounds for spectral gaps. This topic will be considered in a separate paper.\\
Nevertheless the bound is of practical value, because the bounding process $\Z'$ is reversible with respect to $\tilde \pi$, which can be seen by checking the local balance equations.
As a consequence, the bounding techniques for reversible processes, e.g., using Cheeger constants,
found in the literature can be applied directly.\\
In \cite{lorek;szekli:14}[Example 6.2] it is shown, that the bounds obtained via Proposition \ref{prop:gap-bd-abstract}
can be very good for networks with reliable nodes. They compare the bound for an example
provided in \cite{ignatiouk-robert;tibi:12}:\\
This is the network described in Section \ref{sect:networksJackson} with state independent service rates
$\mu_j=\mu$, and routing matrix which fulfills $r_{ii}=0, \forall i=0,1,\dots,J,$
and  $r_{0i}>0, \forall i=1,\dots,J.$
Furthermore, for all $i,j=1,\dots,J,i\neq j,$ holds complete symmetry by $r_{ij}=p\in(0,1)/(j-1)$,
which results in $r_{i0}=1-p(J-1)>0, \forall i=1,\dots,J.$\\
It is assumed that no breakdowns (and repair) occur.\\
In this symmetric network the partial balance (for $D=\emptyset$) holds if and only if
$r_{0i}=1/J, \forall i=1,\dots,J,$ which implies $\lambda_j=\lambda/J, \forall j=1,\dots,J$
and $\eta_j=\lambda/(J(1-p(J-1))), \forall j=1,\dots,J.$\\
We denote by $\mu_{i_0}:=\min\{\mu_{j}:1\leq j\leq J\},$ and,
recalling the bound of the
spectral gap obtained for birth-death processes by van Doorn \cite{doorn:02}, we obtain
from the companion result of Proposition \ref{prop:gap-bd-abstract}
(see  \cite{daduna;szekli:08}[Proposition 4.4])
\begin{equation*}
    Gap(\Z') \geq \left(\sqrt{\mu_{i_0}(1-p(J-1))} - \sqrt{\frac{\lambda}{J}}\right)^2\,.
\end{equation*}
It is easy to check, that for this setting the Assumptions (3.10) and (3.11) of
Corollary 3.4 in \cite{ignatiouk-robert;tibi:12} are fulfilled, which results in an upper bound for $L^2$ spectral gap
\begin{equation*}
    Gap(\Z) \le  \frac{1+p}{1-p(J-2)}\left(\sqrt{\mu_{i_0}(1-p(J-1))} - \sqrt{\frac{\lambda}{J}}\right)^2\,.
\end{equation*}
$\frac{1+p}{1-p(J-2)}$ tends monotonously to $J$ for $p\to 1/(J-1)$, while
for $p\to 0$ it decreases monotonously  to $1$.\\
From the ordering implication {\em
Peskun yields positive definiteness}  follows,
 that if we perturb routing of customers in the networks  by shifting transition probability mass from non diagonal entries into the diagonal (leaving the routing equilibrium fixed)  then the speed of convergence of the perturbed process is smaller.\\
The existence of $L^2$ spectral gap (that is  the question when $Gap(Q^{\Z})>0$) for unreliable networks is a related topic. It is a common knowledge that for networks with constant service rates (not depending on the number of customers at node) the spectral gap for classical Jackson network exists. For service rates that can depend on the number of customers the problem is more delicate.  An {\em iff}  characterization in terms of properties of service rates is given in Lorek and Szekli \cite{lorek;szekli:14}. A special feature of such processes is that the existence of  $L^2$ spectral gap is directly related to the tail properties of the stationary distribution. For references and details see \cite{lorek;szekli:14}.\\
An analogue of Peskun ordering and positive semidefinite order  for generator
matrices is as follows.
\begin{defn}\label{defn:peskun-generator}
Let $Q = (q(x,y):x,y\in \mathbb E)$ and $Q' = (q'(x,y):x,y\in \mathbb E)$ be generator
matrices on a finite set $\mathbb E$ such that $\hat\pi Q=\hat\pi Q'= 0$ holds for a
probability vector $\hat\pi$.\\
$Q'$  is smaller than $Q$ in the positive semidefinite order for generators,
$Q ' \prec_{pd}  Q$, if   $Q- Q' $ is positive semidefinite.\\
$ Q'$  is smaller than $Q$ in the Peskun order for generators,  $ Q' \prec_P Q$,
if for all $x, y\in \mathbb E$ with $x\neq y$  holds $ q'(x,y)\leq q(x,y)$.
\end{defn}
\blemma\label{lem:peskuntopd}
Let $Q = (q(x,y):x,y\in \mathbb E)$ and $Q' = (q'(x,y):x,y\in \mathbb E)$ be generator
matrices on a finite set $\mathbb E$ such that $\hat\pi Q=\hat\pi Q'= 0$ holds for a
probability vector $\hat\pi$. Then
$Q \prec_P Q' \Longrightarrow Q ' \prec_{pd}  Q$ holds.
\bproof
From $ q(x,y)\leq q'(x,y)~\forall x, y\in \mathbb E$ with $x\neq y$ follows
for all $x\in \mathbb E$ that $ q'(x,x)\leq q(x,x)$ holds. So
$Q'-Q := (q'(x,y) - q(x,y):x,y\in \mathbb E)$
is a generator matrix as well. Therefore $-(Q'-Q)$ is positive semidefinite.
\eproof
\elemma
A direct consequence of Definition \ref{defn:peskun-generator}, this lemma, and of Theorem \ref{thm:corrbdr} follows from  \eqref{eq:corrbdr-1}.
\bcor\label{cor:gap-bdr}
Suppose ${\Z}=(\Y,\X)$ and $\Z'=(\Y',\X')$ are ergodic unreliable Jackson network processes,
having the same arrival and service intensities, and with the same routing regimes, described
by  $R=[ r_{ij}]_{i,j=0,1,\dots,J}$ and rerouting fulfilling Assumption \ref{defn:reroutingschemes}.\\
The breakdown-repair process  for $\Z$ is given in Definition \ref{defn:queueindependentbreakdown} and
for $\Z'$ is defined similarly via functions
$A',~B': 2^{\tilde J} \to [0,\infty),$ as given in Theorem \ref{thm:corrbdr}.\\
The processes $\Y$ and $\Y'$ are Markov  with generators $Q^{\Y}=(q^{\Y}(K,L): K, L \subseteq \tilde J)$ of $\Y$
and $Q^{\Y'}=(q^{\Y'}(K,L): K, L \subseteq \tilde J)$ of $\Y'$ as defined in \eqref{eq:bdr-generator} for $\Y$ and similar for $\Y'$.\\
Assume that the stationary distributions of $\Y$ and $\Y'$ are identical, denoted by $\hat\pi$.\\
If $Q^{\Y} \prec_{pd}  Q^{\Y'}$ holds, then
$Gap(Q^{\Y'}) \le  Gap(Q^{\Y}).$
\ecor
An easy to understand property is that whenever the breakdown-repair process $\Y$ of $\Z$
is uniformly faster than the breakdown-repair process $\Y'$ of $\Z'$,
i.e., for all
$x\neq y$  holds $ q^{\Y'}(x,y)\leq q^{\Y}(x,y)$, we have
$Gap(Q^{\Y'})  \leq Gap(Q^{\Y}).$ This follows directly from Lemma \ref{lem:peskuntopd}.\\
So, for example, if we have ${A(D)} = \kappa^{|D|}\cdot A'(D),$
and ${B(D)} = \kappa^{|D|}\cdot B'(D), A\in 2^{\tilde J},$ for some $\kappa > 1$,
then $Q^{\Y'} \prec_{P}  Q^{\Y}$, and these breakdown-repair processes fulfill the requirement
of Corollary \ref{cor:gap-bdr}.
\subsection{Asymptotic variance}\label{sect:asymptvariance}
Peskun \cite{peskun:73} and Tierney \cite{tierney:98} derived comparison theorems for the asymptotic variance of  Markov chains for application to optimal selection of MCMC transition kernels in discrete time.
These asymptotic variances occur as variance in the limiting distribution of central limit theorems (CLTs)
for the MCMC estimators. For our network processes $\Z$ we consider  Markov chain $(X_k,k\ge 1)$, say
with transition matrix $K= Id + \varepsilon Q^\Z$
(with $\varepsilon>0$ sufficiently small). Under some regularity conditions  on a homogeneous Markov chain with one step transition kernel
$K$ we can obtain CLT  of the form
$$
\sqrt {n} (\frac{1}{n} \sum_{k=1}^n f(X_k) - E_{\tilde\pi}(f(X_t))) \stackrel{D}{\to} N(0,v(f,K)),
$$
where the asymptotic variance is
$v(f,K) = \langle f,f\rangle_{\tilde\pi} -\tilde\pi(f) + 2 \sum_{k=1}^\infty \langle
f,K^kf\rangle_{\tilde\pi}.$\\
Regularity conditions under which CLT holds for such Markov chains is a  topic which we shall study in a separate paper. For reversible chains with positive spectral gap it is possible to give conditions in terms of the service rates,  but a general {\em iff} characterization in terms of the service rates seems to be an open  problem.
\begin{prop}\label{prop:varorder-nets}
 Consider two ergodic unreliable Jackson networks with the same arrival and service
intensities, and  state processes  $\Z$ and $\Z'$. Assume
that the  routing matrices $R$ and $R'$ are reversible
with respect to $\xi$.
Both networks follow a  rerouting mechanism according to \eqref{eq:routing-rejectionD2}
with the property \eqref{eq:routing-rejectionD3}, such that $R^D$ and $R'^D$ are reversible
with respect to $\xi^D$.\\
If  $R^D$ and $R'^D$ are ordered for all $D$ in  positive definite  order, $R'^D \prec_{pd}R^D$,
then for any  function $f\in L^2_0(\tilde{\mathbb E},\tilde \pi) :=
\{g \in L^2(\tilde{\mathbb E},\tilde \pi):\pi(g)=0\}$ holds
$v(f,Id +\varepsilon Q^{\Z}) \geq v(f,Id +\varepsilon Q^{{\Z}'}).$
\end{prop}
\bproof
For standard Jackson networks without breakdown and repair it is well known that reversibility
 of the routing matrix $R$ implies reversibility of the joint queue length process.
 A direct way to prove this is to check the local balance equations with respect to
 the stationary distribution $\pi$. It is easy to see that this way of proof verifies reversibility
 of the processes $\Z$ and $\Z'$ here as well. The reason is that the breakdown and repair process
 is reversible with respect to $\hat \pi$, and that for fixed $D$ and $\ourN_D$ intensities of
 possible transitions  on $\ourN^{|\tilde J \setminus D|}$ balance locally with respect to
the densities
$\prod_{\ell\in \tilde J\setminus D} \prod_{i=1}^{n_\ell} \left(\frac{\eta_\ell}{\mu_\ell(i)}\right)$.\\
Because $\Z$ and $\Z'$ are irreducible we can apply a result of Mira and Geyer
\cite{mira;geyer:99}[Theorem 4.2], which states that under this condition
the required inequality is equivalent to ordering of the one step correlations for $f\in L^2_0(\mathbb E,\tilde \pi)$.
The letter statement can be shown exactly as in the proof of Proposition \ref{prop:pdorder-nets}.
\eproof

\section{Jackson networks with finite waiting rooms}\label{sect:finite-buffers}
We consider the Jackson networks from Section \ref{sect:networksJackson} and assume now
that some or all of
the single server nodes with queue length dependent service intensities have a
finite capacity.
This means that such nodes have a service place and only a finite number of waiting places available to store waiting customers. 
The problem one is faced with is that  whenever a node, say $j$, reaches its maximal buffer size
then no further arriving customer can enter node $j$ and we have to ''resolve'' such blocking situation.\\
There are many ''blocking protocols'' available in practice and in theoretical models but the standard
observation is that when reducing the infinite buffers of the Jackson networks from Section \ref{sect:networksJackson} to finite sizes the stationary distribution is no longer of product form,
in fact in almost all cases  the stationary distribution is not available. The  simplest reaction
to a full buffer situation with a newly arriving customer is to discard this customer - this leads
to a simple equilibrium only in case of a single station network.\\
A survey, emphasizing applications in telecommunication systems is given in
\cite{perros:90}[Section 2], a general overview is \cite{balsamo;denittopersone;onvural:01}.\\
It turned out that for obtaining product form stationary distributions for finite buffer Jackson networks
one has to introduce clever rerouting strategies, i.e. to redirect and redistribute newly arriving customers who want to enter nodes with full buffers in such a way that the the system's customer flow is
balanced.
A literature review shows that the most prominent rerouting strategies in case of full buffers
which lead to product form equilibrium are the
schemes mentioned in Example \ref{ex:rerouting}. These are described in more detail
e.g. in \cite{dijk:11}[Section 1.6.3]
under the headings {\sc Conservative protocol} $\equiv$ Stalling, {\sc Jump-over protocol} $\equiv$
Skipping, and in \cite{perros:90}[p. 456,457] as {\sc repetitive-service-random destination}.
It turns out that in all three cases the stationary distribution for feasible states in the finite buffer
network are (up to normalization) of the form \eqref{eq:steadystateJ}.
A detailed study of the jump-over protocol in case of blocking is \cite{dijk:88}.\\
As a referee pointed out, it is therefore a tempting conjecture that correlation formulas similar to ours should be valid in the context of finite buffer networks as well. But going into the details shows
that there are subtleties in the structure and the functioning of the finite buffer networks which
must be handled careful. 

{\bf (1)} The network of Section \ref{sect:networksJ-unrel} is influenced by an external environment
which is Markov for its own, and changes of the environment enforces the Jackson network to redirect
 its customers, if necessary. In case of the finite buffer network the customers are redirected by
 intrinsic forces (state space restrictions).\\
{\bf (2)} The intrinsic forces which influence routing and rerouting  lead to queue
length dependent routing decisions. Seemingly, this makes it nearly impossible to formulate
an overall simple condition in parallel to Assumption \ref{defn:reroutingschemes}, which produces
product form stationary distributions.\\
{\bf (3)} The area of applications (problems)  which is the classical source of interest in finite buffer
networks with blocking protocols is in telecommunications networks. There it is unrealistic to assume that
in case of blocking (=full buffer) the service at the full buffer node is interrupted, which is natural
at broken down nodes in the setting of Section \ref{sect:networksJ-unrel}.

Nevertheless, it turned out that correlation formulas similar to those in Propositions
\ref{prop:corr-formula} and \ref{prop:corrformula-split} can be derived in case of the most popular
unblocking scheme  (repetitive-service-random destination) for networks with reversible routing.
This and more general results are part of our ongoing research and will be summarized in a
forthcoming technical report \cite{daduna;szekli:15}.

\section{Proof of the correlation formulas}\label{sect:proofs}
For $f,g: \ourEE\to \R$ and steady state probability $\tilde\pi$ of ${\Z}$
we are interested in the one-step correlation expressions
\beqar\label{eq:corr-formula-0-pr1}
&&\langle f,Q^\Z  g\rangle_{\tilde\pi}=C^{-1}\sum_{D\subseteq {\{1,\dots,J\}}} \sum_{\n\in \ourN^J}
\hat\pi(D)
\prod_{\ell=1}^J \prod_{i=1}^{n_\ell} \left(\frac{\eta_\ell}{\mu_\ell(i)}\right) f(D,\n)
 \left\{\sum_{j\in \tilde J \setminus D} \lambda^D_j g(D,\n+e_j)\right.\\
&&+ \sum_{j\in \tilde J \setminus D} \mu_j(n_j) r^D_{j0} g(D,\n-e_j) 
+ \sum_{j\in \tilde J \setminus D}
\sum_{{i\in \tilde J \setminus D},{i\neq j}}\mu_j(n_j) r^D_{ji} g(D,\n-e_j+e_i)\nonumber\\
&&-\left(\sum_{j\in \tilde J \setminus D} \lambda^D_j
+ \sum_{j\in \tilde J \setminus D} \mu_j(n_j)(1-r^D_{jj})\right)g(D,\n)\nonumber\\
&&+\sum_{H\subset D} q^{\Y}(D, H) g(H,\n)
+ \sum_{I\supset D} q^{\Y}(D, I)g( I,\n)
\left.-\left(\sum_{H\subset D} q^{\Y}(D, H)
+ \sum_{I\supset D} q^{\Y}(D,  I)
\right)g(D,\n)
\right\}.\nonumber
\eeqar
\bproof (of Proposition \ref{prop:corrformula-split})
Interchanging summations, regrouping terms, and exploiting the
product form structure of the state distribution, we obtain from
\eqref{eq:corr-formula-0-pr1}
\beqaro
&&\sum_{\n\in \ourN^J} \pi(\n) \left[\sum_{D\subseteq {\{1,\dots,J\}}} \hat\pi(D)\ \   
 \left\{\sum_{H\subset D} q^{\Y}(D,H) f(D,\n) g(H,\n)
 + \sum_{I\supset D} q^{\Y}(D, I) f(D,\n) g(I,\n)
\right\}\right.\nonumber\\
&&\qquad\qquad\qquad\qquad\left.- f(D,\n) g(D,\n)
\left\{\sum_{H\subset D} q^{\Y}(D, H)  + \sum_{I\supset D} q^{\Y}(D, I)\right\}\right]\\
&&+ \sum_{D\subseteq {\{1,\dots,J\}}} \hat\pi(D)
\sum_{\n_D\in \ourN^{|D|}}
\prod_{\ell\in D} \left( C_\ell^{-1}\prod_{i=1}^{n_\ell} \frac{\eta_\ell}{\mu_\ell(i)}\right)
\eeqaro
\beqaro
&&\left[
\sum_{\n_{\tilde J\setminus D}\in \ourN^{J-|D|}}
\prod_{\ell\in \tilde J\setminus D} \left( C_\ell^{-1}
\prod_{i=1}^{n_\ell} \frac{\eta_\ell}{\mu_\ell(i)}\right)
f(D,(\n_D,\n_{\tilde J \setminus D}))\left\{\sum_{j\in \tilde J \setminus D} \lambda^D_j g(D,\n+e_j)\right.\right.\\
&&
+ \sum_{j\in \tilde J \setminus D} \mu_j(n_j) r^D_{j0}
g(D,(\n_D,\n_{\tilde J \setminus D}-e_j))
+ \sum_{j\in \tilde J \setminus D}
\sum_{{i\in \tilde J \setminus D},{i\neq j}}\mu_j(n_j) r^D_{ji}
g(D,(\n_D,\n_{\tilde J \setminus D}-e_j+e_i))\nonumber\\
&&\qquad\qquad\qquad\qquad-\left.\left.\left(\sum_{j\in \tilde J \setminus D} \lambda^D_j
+ \sum_{j\in \tilde J \setminus D} \mu_j(n_j)(1-r^D_{jj})\right)
g(D,(\n_D,\n_{\tilde J \setminus D}))
\right\}\right]
\eeqaro
For each fixed $D,\n_D,$ the terms in the last squared brackets are identical to the one step correlation
of a Jackson network in equilibrium on node  set $\tilde J \setminus D$
(with the respective transition rates) with respect to the functions
$f(D,(\n_D,{(\cdot)}_{\tilde J \setminus D}))$ and $g(D,(\n_D,{(\cdot)}_{\tilde J \setminus D}))$.\\
We have agreed to denote the generator of such network by $Q^{\X_{\tilde J \setminus D}}$,
and its steady state by $\pi_{\tilde J \setminus D}$, which leads to the proposed formula with
the aid of the synthetic networks.
\eproof
\bproof (of Proposition \ref{prop:corr-formula})
We restart with the expression \eqref{eq:corr-formula-0-pr1} and
observe that for fixed ${D\subseteq {\{1,\dots,J\}}}$ the contribution of $-r^D_{jj}$ in the
negative terms would be exactly the contribution in the double sum of
${i\in \tilde J \setminus D},{i = j}$ in the positive terms, where for $i = j$ would occur $g(D,\n-e_j+e_j)=g(D,\n)$ otherwise. Together with $\mu_j(0)=0 \  \forall j$, incorporating these contributions simplifies our expression to
\beqar\label{eq:corrnet2}
&& -C^{-1}\sum_{D\subseteq {\{1,\dots,J\}}} \sum_{\n\in \ourN^J} \hat\pi(D)
\prod_{\ell=1}^J \prod_{i=1}^{n_\ell} \left(\frac{\eta_\ell}{\mu_\ell(i)}\right) f(D,\n) g(D,\n)\\
&&\qquad\qquad\qquad\left[\left(\sum_{H\subset D} q^{\Y}(D,H)+ \sum_{I\supset D} q^{\Y}(D, I)\right)
+\left(\sum_{j\in \tilde J \setminus D} \lambda^D_j
+ \sum_{j\in \tilde J \setminus D}
\mu_j(n_j)\right) \right]\nonumber\\
&&+ C^{-1}\sum_{D\subseteq {\{1,\dots,J\}}} \sum_{\n\in \ourN^J}  \hat\pi(D)
\prod_{\ell=1}^J \prod_{i=1}^{n_\ell} \left(\frac{\eta_\ell}{\mu_\ell(i)}\right)
f(D,\n)
\left\{\sum_{H\subset D} q^{\Y}(D, H) g( H,\n) +
\sum_{I\supset D} q^{\Y}(D, I) g( I,\n)
\right\}\nonumber\\
&&+ C^{-1}\sum_{D\subseteq {\{1,\dots,J\}}} \sum_{\n\in \ourN^J} \hat \pi(D)
\prod_{\ell=1}^J \prod_{i=1}^{n_\ell} \left(\frac{\eta_\ell}{\mu_\ell(i)}\right)
f(D,\n)
\left\{\sum_{j\in \tilde J \setminus D} \lambda^D_j g(D,\n+e_j)\right.\nonumber\\
&&\qquad\qquad\qquad+ \sum_{j\in \tilde J \setminus D}
\mu_j(n_j) r^D_{j0} g(D,\n-e_j) 
\left.+ \sum_{j\in \tilde J \setminus D}
\sum_{{i\in \tilde J \setminus D}}\mu_j(n_j) r^D_{ji} g(D,\n-e_j+e_i)\right\}\nonumber
\eeqar
The last two lines of the formula above turn into
\beqaro
&& C^{-1}\sum_{D\subseteq {\{1,\dots,J\}}} \hat\pi(D) \sum_{\n_D\in \ourN^D}
\prod_{\ell\in D} \prod_{i=1}^{n_\ell} \left(\frac{\eta_\ell}{\mu_\ell(i)}\right)\\
&&\left[
\sum_{\n_{\tilde J \setminus D}\in \ourN^{\tilde J \setminus D}}
\prod_{\ell\in {\tilde J \setminus D}} \prod_{i=1}^{n_\ell} \left(\frac{\eta_\ell}{\mu_\ell(i)}\right)
f(D,\n)
\left\{\sum_{j\in \tilde J \setminus D} \lambda^D_j g(D,\n+e_j)\right.\right.\\
&&\qquad\qquad\qquad\left.\left.+ \sum_{j\in \tilde J \setminus D} 
\mu_j(n_j) r^D_{j0} g(D,\n-e_j)
+ \sum_{j\in \tilde J \setminus D} 
\sum_{{i\in \tilde J \setminus D}}\mu_j(n_j) r^D_{ji} g(D,\n-e_j+e_i)\right\}
\right.
\eeqaro
\beqaro
&&= C^{-1}\sum_{D\subseteq {\{1,\dots,J\}}} \hat\pi(D) \sum_{\n_D\in \ourN^D}
\prod_{\ell\in D} \prod_{i=1}^{n_\ell} \left(\frac{\eta_\ell}{\mu_\ell(i)}\right)\\
&&\left[
\sum_{\n_{\tilde J \setminus D}\in \ourN^{\tilde J \setminus D}}
\prod_{\ell\in {\tilde J \setminus D}} \prod_{i=1}^{n_\ell} \left(\frac{\eta_\ell}{\mu_\ell(i)}\right)
\left\{\sum_{j\in \tilde J \setminus D} f(D,\n) \lambda g(D,\n+e_j) {r^D_{0j}}\right.\right.\\
&&\left.\left.+ \sum_{j\in \tilde J \setminus D} 
\mu_j(n_j)  f(D,\n) g(D,\n-e_j) r^D_{j0}
+ \sum_{j\in \tilde J \setminus D} 
\sum_{{i\in \tilde J \setminus D}}\mu_j(n_j) f(D,\n) g(D,\n-e_j+e_i) r^D_{ji}\right\}
\right]
\eeqaro
In the last line, if $n_j>0$, the expression $\mu_j(n_j)$ cancels against a factor in the steady state probability. This leads to
(we underbrace some intuition and use $\eta_j=\eta_j^D$ for $j\in \tilde J\setminus D$)
\beqar\label{eq:rej-1}
&& \sum_{D\subseteq {\{1,\dots,J\}}} \sum_{\n\in \ourN^J} \tilde\pi(D,n)
\left\{\sum_{j\in \tilde J \setminus D}  \underbrace{\lambda f(D,\n) g(D,\n+e_j)
{r^D_{0j}}}_{0\to j}
+ \sum_{j\in \tilde J \setminus D} 
 \underbrace{\eta_j^D  f(D,\n+e_j) g(D,\n) r^D_{j0}}_{j\to 0}\right.\\
&&\left.+ \sum_{j\in \tilde J \setminus D} 
\sum_{{i\in \tilde J \setminus D}} \underbrace{\eta_j^D f(D,\n+e_j) g(D,\n+e_i)
r^D_{ji}}_{j\to i}
+  \underbrace{\lambda f(D,\n) g(D,\n) {r^D_{00}}}_{0\to 0} -\lambda f(D,\n) g(D,\n) {r^D_{00}}
\right\}\nonumber
\eeqar
From Lemma \ref{lem:reroutingschemesR},  with $\eta^D_{0}:= \lambda$ the vector
$\hat\eta^D :=(\eta^D_j, j\in \tilde J_0\setminus D)$  solves  $x=x\cdot R^D$.
Insert the normalized solution $\xi^D := (\xi^D_i:i\in  \tilde J_0\setminus D)$
 into \eqref{eq:rej-1}, and
then insert into the  correlation expressions
\eqref{eq:corrnet2} to obtain
\begin{eqnarray*}
&&\langle f,Q^\Z g\rangle_{\tilde\pi}
 = \sum_{D\subseteq {\{1,\dots,J\}}} \sum_{\n\in \ourN^J} \tilde\pi(D,\n)
 \left\{\sum_{H\subset D} q^{\Y}(D, H) f(D,\n) g( H,\n)
 + \sum_{I\supset D} q^{\Y}(D, I) f(D,\n) g( I,\n)
\right\}\nonumber\\
&&+ \sum_{D\subseteq {\{1,\dots,J\}}} \sum_{\n\in \ourN^J} \tilde\pi(D,\n) \frac{\lambda}{\xi^D_0}
\left\{\sum_{j\in (\tilde J \cup \{0\}) \setminus D} \sum_{i\in (\tilde J \cup \{0\}) \setminus D}
\xi^D_j r^D_{ji} f(D,\n+e_j) g(D,\n+e_i)\right\} \nonumber\\
&&- \sum_{D\subseteq {\{1,\dots,J\}}} \sum_{\n\in \ourN^J} \tilde\pi(D,\n)
\lambda f(D,\n) g(D,\n) {r^D_{00}}\nonumber\\
&&- \sum_{D\subseteq {\tilde J}} \sum_{\n\in \ourN^J} \tilde\pi(D,\n) f(D,\n) g(D,\n)
\left\{\left[\sum_{H\subset D} q^{\Y}(D, H)  + \sum_{I\supset D} q^{\Y}(D, I) \right]\right.
\left.+\left[\sum_{j\in \tilde J \setminus D} \lambda_j^D + \sum_{j\in \tilde J \setminus D}
\mu_j(n_j) \right]\right\}\nonumber
\end{eqnarray*}
With $\sum_{j\in \tilde J \setminus D} \lambda^D_j +\lambda r^D_{00}{=\lambda}$
this yields finally the desired correlation formula from Proposition \ref{prop:corr-formula}.
\eproof

{\bf Acknowledgement.} We thank an anonymous referee for careful reading of the first version of the paper and her or
his constructive critics.

\end{document}